\documentclass{amsart}
\usepackage{amsthm,amssymb}

\newtheorem{theorem}{Theorem}

\newtheorem{lemma}[theorem]{Lemma}

\newtheorem{remark}{Remark}

\newtheorem{conjecture}{Conjecture}

\begin{document}
\title[Polynomials ]{On Erd\'{e}lyi-Magnus-Nevai conjecture for Jacobi polynomials}
\author[I. Krasikov]{Ilia Krasikov}
\address{ Department of Mathematical Sciences, Brunel University, Uxbridge
UB8 3PH United Kingdom} \email{mastiik@brunel.ac.uk}

\begin{abstract}
T. Erd\'{e}lyi, A.P. Magnus and P. Nevai conjectured that for
$\alpha , \beta \ge - \, \frac{1}{2} \, ,$  the orthonormal Jacobi
polynomials ${\bf P}_k^{( \alpha , \beta )} (x)$ satisfy the
inequality
\begin{equation*}
\max_{x \in
[-1,1]}(1-x)^{\alpha+\frac{1}{2}}(1+x)^{\beta+\frac{1}{2}}\left({\bf
P}_k^{( \alpha , \beta )} (x) \right)^2 =O \left(\max
\left\{1,(\alpha^2+\beta^2 )^{1/4} \right\} \right),
\end{equation*}
[Erd\'{e}lyi et al.,Generalized Jacobi weights, Christoffel
functions, and Jacobi polynomials, SIAM J. Math. Anal. 25 (1994),
602-614]. Here we will confirm this conjecture in the
ultraspherical case $\alpha = \beta \ge \frac{1+ \sqrt{2}}{4},$
even in a stronger form by giving very explicit upper bounds. We
also show that
\begin{equation*}
\sqrt{\delta^2-x^2} \, (1-x^2)^{\alpha}\left({\bf P}_{2k}^{(
\alpha , \alpha)} (x)\right)^2 < \frac{2}{\pi} \left( 1+
\frac{1}{8(2k+ \alpha)^2}  \right)
\end{equation*}
for a certain choice of $\delta,$ such that the interval $(-
\delta, \delta)$  contains all the zeros of ${\bf P}_{2k}^{(
\alpha , \alpha)} (x).$ Slightly weaker bounds are given for
polynomials of odd degree.

\hspace{1pt}

\end{abstract}

\maketitle


\section{Introduction}
\label{secintr}
 In this paper we will use bold letters for
orthonormal polynomials versus  regular characters for orthogonal
polynomials in the standard normalization \cite{szego}.

 Given a family $\{{\bf
p}_i(x) \}$ of orthonormal polynomials orthogonal on a finite or
infinite interval $I$ with respect to a weight function $w(x) \ge
0,$ it is an important and difficult problem to estimate $\sup_{x
\in I}\sqrt{w(x)} \, |{\bf p}_i(x)|,$ or, more generally, to find
an envelope of the function $\sqrt{w(x)} \,{\bf p}_i(x)$ on $I.$
Those two questions become almost identical if we introduce an
auxiliary function $\phi(x)$ such that $\sqrt{\phi(x)w(x)} \,{\bf
p}_i(x)$ exhibits nearly equioscillatory behaviour. Of course, the
existence of such a function is far from being obvious but it
turns out that in many cases one can choose $\phi =\sqrt{(x-
d_m)(d_M-x)} \, ,$ with $d_m, d_M$ being appropriate
approximations to the least and the largest zero of $p_i$
respectively. The simplest example is given by Chebyshev
polynomials $ T_i (x)$ and $\phi=\sqrt{1-x^2}.$ This illustrates a
classical result of G.\,Szeg\"{o}  asserting that for a vast class
of weights on $[-1,1]$ and $i \rightarrow \infty ,$ the function
$\sqrt{\sqrt{1-x^2} \, w(x)} \, {\bf p}_i(x)$ equioscillates
between $\pm \, \sqrt{\frac{2}{\pi}} \, ,$ \cite{szego}.

A very general theory for exponential weights $w=e^{-Q(x)}$
stating that under some technical conditions on $Q,$
$$
\max_{I} \left| \sqrt{\sqrt{|(x- a_{-i})(a_i-x)|} \, w(x)} \, {\bf
p}_i(x) \right|< C,
$$
where the constant $C$ is independent on $i$ and $a_{\pm i}$ are
Mhaskar-Rahmanov-Saff numbers for $Q,$ was developed by A.L. Levin
and D.S. Lubinsky \cite{levlub}. Recently it has been extended to
the Laguerre-type exponential weights $x^{2 \rho}e^{-2 Q(x)}$
\cite{kasuga,levin4}.

In the case of classical orthogonal Hermite and Laguerre
polynomials explicit bounds confirming such a nearly
equioscillatory behaviour independently on the parameters involved
were given in \cite{krasherm} and \cite{krasl2} respectively.

The case of Jacobi polynomials $P_k^{(\alpha , \beta )} (x), \;
\;w(x)=(1-x)^{\alpha} (1+x)^{\beta } ,$ is much more difficult.
Let us introduce some necessary notation.

We define
$$M_k^{ \alpha , \beta}(x,d_m , d_M ) =
 \sqrt{(x-d_m)(d_M-x)} \; (1-x)^{ \alpha }(1+x)^{\beta} \left( {\bf P}_k^{(\alpha
, \beta )} (x)\right)^2 ,$$
 $$\mathcal{M}_k^{ \alpha , \beta}(d_m , d_M ) =
 \max_{x \in [-1,1]} M_k^{ \alpha , \beta}(x;d_m , d_M ) ,$$
 what we will abbreviate to $M_k^{ \alpha , \beta}(x)$ and $\mathcal{M}_k^{ \alpha , \beta}$ if
 $d_m=-1, \; d_M=1,$ that is for $\phi(x)=\sqrt{1-x^2}
 \, .$ We will also omit one of the superscripts in the
 ultraspherical case $\alpha = \beta$ writing, for example, $M_k^{ \alpha
 }(x)$ instead of $M_k^{ \alpha , \alpha}(x),$ and shorten $M_k^{ \alpha
 }(x, -d , d), \; \mathcal{M}_k^{ \alpha
 }( -d , d)$ to $M_k^{ \alpha
 }(x,  d), \; \mathcal{M}_k^{ \alpha
 }( d)$ respectively.

 As $P_k^{(\alpha, \beta )} (x)=(-1)^k P_k^{(\beta , \alpha)}
 (-x)$ we may safely assume that $ \alpha \ge \beta.$

 For $- \,\frac{1}{2} <\beta \le \alpha < \frac{1}{2} \, ,$ the following is known \cite{chow}:
\begin{equation}
\label{gatwong} \mathcal{M}_k^{ \alpha , \beta}  \le \frac{2^{2
\alpha +1}\Gamma (k+ \alpha+ \beta +1) \Gamma (k+ \alpha +1)}{\pi
 k! \,
 (2k+ \alpha +\beta +1)^{2 \alpha} \Gamma (k+ \beta+1)}= \frac{2}{\pi}+O\left(\frac{1}{k} \right)\, ,
 \end{equation}
where $k=0,1,...\,$.

A slightly stronger inequality in the ultraspherical case was
obtained earlier by L. Lorch \cite{lorch}.

 A remarkable result covering almost all possible range of the parameters
 has been established by T.\,Erd\'{e}lyi, A.P.\,Magnus and P.\,Nevai,
\cite{erdel},

\begin{equation}
\label{eqernt}
 \mathcal{M}_k^{ \alpha , \beta}  \le
 \frac{2e \left(2+ \sqrt{\alpha^2+\beta^2} \right)}{\pi} \, ,
\end{equation}
provided $k \ge 0, \; \; \alpha, \beta \ge - \frac{1}{2} \, .$

Moreover, they suggested the following conjecture:
 \begin{conjecture}
 \label{conjemn}
 $$\mathcal{M}_k^{ \alpha , \beta} = O \left( \max \left \{1, | \alpha |^{1/2} \right\} \,
\right),$$ provided $ \alpha \ge \beta \ge - \, \frac{1}{2} \, .$
\end{conjecture}

The best currently known bound was given by the author
\cite{krasjac},
\begin{equation}
\label{krjac} \mathcal{M}_k^{ \alpha , \beta}  \le 11 \, \left(
\frac{(\alpha+ \beta+1)^2(2k+\alpha+ \beta+1)^2}{4k(k+ \alpha+
\beta+1)} \right)^{1/3}=O \left( \alpha^{2/3} \left(1+
\frac{\alpha}{k} \right)^{1/3} \right),
\end{equation}
provided $k \ge 6, \; \; \alpha \ge \beta  \ge \frac{1+
\sqrt{2}}{4} \, .$

We also brought some evidences in support of the following
stronger conjecture
\begin{conjecture}
 \label{conjik}
 $$\mathcal{M}_k^{ \alpha , \beta} =O \left( \max \left \{1, |\alpha |^{1/3} \left(1+
\frac{|\alpha |}{k} \right)^{1/6} \right\} \, \right),$$ provided
$ \alpha \ge \beta  \ge - \, \frac{1}{2} \,.$
\end{conjecture}
Here we will confirm this conjecture in the ultraspherical case.
Namely we prove the following
\begin{theorem}
\label{mainth1} Suppose that $k \ge 6, \; \; \alpha = \beta  \ge
\frac{1+ \sqrt{2}}{4} \, .$ Then
\begin{equation}
\label{gleq1}
 \mathcal{M}_k^{ \alpha } <\mu \, \alpha ^{1/3}
\left(1+ \frac{\alpha }{k} \right)^{1/6},
\end{equation}
where
$$
\mu = \left\{
\begin{array}{cc}
\frac{10}{7} \, , & k \; \; even,\\
& \\
22, & k \; \; odd.
\end{array}
\right.
$$
\end{theorem}
We deduce this result from the following two theorems. The first,
which has been established in \cite{krasjac}, gives a sharp
inequality for the interval containing all the local maxima of the
function $M_k^{ \alpha , \beta}(x).$  The second one will be
proven here and in fact demonstrates equioscillatory behaviour of
$M_k^{\alpha } (x, d )$ under an appropriate choice of $d.$
\begin{theorem}
\label{grmax} Suppose that $k \ge 6, \; \; \alpha \ge \beta \ge
\frac{1+ \sqrt{2} }{4} \, .$ Let $x$ be a point of a local
extremum of $
 M_k^{ \alpha, \beta} (x).$ Then $x \in \left(\eta_{-1}, \eta_1
 \right),$
 where
 \begin{equation}
 \label{eqN}
 \eta_j= j \left(\cos (\tau+j \omega)-
 \theta_j \left( \frac{\sin^4 (\tau+j \omega)}{2 \cos \tau \cos \omega} \right)^{1/3}
 (2k+ \alpha+ \beta+1)^{-2/3} \right)
 \end{equation}
$$
\sin \tau= \frac{\alpha+\beta+1}{2k+ \alpha+ \beta+1} \, , \; \;
\sin \omega= \frac{\alpha- \beta}{2k+\alpha+ \beta+1} \, , \; \; 0
\le \tau , \, \omega < \frac{\pi}{2} \, ;
$$
 and
$$
\theta_j = \left\{
\begin{array}{cc}
1/3 \, , & j=-1,\\
&\\
3/10 \, , & j=1.
\end{array}
\right.
$$
 In particular, in the ultraspherical case
\begin{equation}
 \label{eqNus}
|x| < \eta = \cos \tau \left(1 -  \frac{  2^{-1/3}}{3} \, \, (2k+2
\alpha+1)^{-2/3} \, tan^{4/3} \tau
 \right),
\end{equation}
with $\sin \tau= \frac{2 \alpha+1}{2k+2 \alpha+1} \, . $
\end{theorem}

\begin{theorem}
\label{mainth2} Suppose that $\alpha > \frac{1}{2} \, ,$ and let
\begin{equation}
\label{defdelta} \delta= \sqrt{1- \frac{4\alpha^2-1}{(2k+2
\alpha+1)^2-4}} \, .
\end{equation}
Then
\begin{equation}
\label{gleq2}  \mathcal{M}_k^{ \alpha }( \delta ) < \left\{
\begin{array}{cc}
\frac{2}{\pi} \left(1+ \frac{1}{8(k+ \alpha)^2 }\right), & k \ge
2, \; \;  even, \\
&\\
\frac{230}{\pi} \, , & k \ge 3, \; \; odd.
\end{array}
\right.
\end{equation}
Moreover, all local maxima of the function $M_k^{\alpha }(x)$ lie
inside the interval $(- \delta , \delta ).$
\end{theorem}
To prove this theorem we construct an envelope of $M_k^{ \alpha ,
\beta}(x; d_m , d_M )$ using so-called Sonin's function. Then we
show that in the ultraspherical case for $\alpha > \frac{1}{2} $
it has the only minimum at $x=0$ if $\delta_m=-1, \; \delta_M=1, $
whereas for $-d_m = d_M =\delta$ the point $x=0$ is the only
maximum. Sharper bounds for the even case are due to the fact that
$x=0$ is the global maximum of $M_{2k}^{\alpha }(x, \delta )$ and
the value of $P_{2k}^{(\alpha , \alpha )} (0)$ is known.

The paper is organized as follows. In the next section we present
a simple lemma being our main technical tool. We will illustrate
it by proving that the function $M_k^{\alpha \, \beta }(x)$ is
unimodal with the only minimum in a point depending only on
$\alpha$ and $\beta .$ The even and the odd cases of Theorem
\ref{mainth2} will be proven in sections \ref{seceven} and
\ref{secodd} respectively. The last section deals with the proof
of Theorem \ref{mainth1}.
\section{Preliminaries}
\label{secprel}
 In his seminal book \cite{szego} Szeg\"{o} presented
a few result concerning the behaviour of local extrema of
classical orthogonal polynomials based on an elementary approach
via so-called Sonin's function. In particular, he gave a
comprehensive treatment of the Laguerre polynomials \cite[Sec
7.31, 7.6 ]{szego}, but did not try to deal with the Jacobi case
for arbitrarily values of $\alpha$ and $\beta.$ Here we combine
his approach with the following very simple idea.

Given a real function $f(x),$ Sonin's function $S=S(f;x)$ is
$S=f^2+ \psi (x){f'^2},$ where $\psi(x) >0$ on an interval
$\mathcal{I}$ containing all local maxima of $f.$
 Thus, they lie on $S,$ and if $S$ is unimodal we can locate the global one.

\begin{lemma}
\label{mainl} Suppose that a function $f$ satisfies on an open
interval $\mathcal{I}$ the Laguerre inequality
 \begin{equation}
 \label{lagineq}
 f'^2-f f'' >0,
 \end{equation}
 and a differential equation
\begin{equation}
\label{difeqsz}
 f''-2A(x)f'+B(x)f=0,
 \end{equation}
  where $A \in \mathbb{C}(\mathcal{I}), \; \; B(x) \in \mathbb{C}^1(\mathcal{I}),$
  and $B$ has at most two zeros on $\mathcal{I} .$
 Let
$$
S(f;x)=f^2+\frac{f'^2}{B},
$$
 then all the local maxima of $f$ in $\mathcal{I}$ are in the
 intervals
 defined by $B(x) >0,$ and
 $$
 Sign \left( \frac{d}{dx} S(f;x) \right) =Sign (4A B-B').
 $$
\end{lemma}
\begin{proof}
We have $0 <f'^2-f f''=f'^2-2A f f'+B f^2,$ hence $B(x) >0$
whenever $f' =0.$ Finally,
$$
\frac{d}{dx} \left(f^2+\frac{f'^2}{B} \right) = \frac{4A
B-B'}{B^2} \, f'^2(x),
$$
and $B(x) \ne 0$ in one or two intervals containing all the
extrema of $f$ on $\mathcal{I}.$
\end{proof}
Let us make a few remarks concerning the Laguerre inequality
(\ref{lagineq}). Usually it is stated for hyperbolic polynomials,
that is real polynomials with only real zeros, and their limiting
case, so-called Polya-Laguerre class. In fact, it holds for a much
vaster class of functions. Let $L(f)=f'^2-f f'',$ defining
$\mathcal{L} =\{f(x):L(f)
>0 \},$ we observe that $\mathcal{L}$ is closed under linear
transformations $x \rightarrow a x+b.$ Moreover, since
$$
L(f g)=f^2 L(g)+g^2 L(f),
$$
$\mathcal{L}$ is closed under multiplication as well. Thus,
$L(x^{\alpha} )= \alpha x^{2 \alpha -2},$ yields the polynomial
case and much more. Many examples may be obtain by $L \left( e^{f}
\right)=-e^{2 f} f''$ and obvious limiting procedures.

For our purposes it is enough that (\ref{lagineq}) holds for the
functions
$$ \left((x-d_m)(d_M-x)
\right)^{1/4}(1-x)^{\alpha/2} (1+x)^{\beta/2}P_k^{( \alpha , \beta
)} (x),$$ provided $-1 \le d_m <x <d_M \le 1,$ and $ \alpha ,
\beta \ge 0.$

 To demonstrate how powerful this lemma is, we apply it to $M_k^{\alpha
, \beta }(x)$ to show that its local maxima lie on a unimodal
curve.

 From the differential equation for Jacobi polynomials
\begin{equation}
\label{difeqy} (1-x^2)y''=((\alpha+\beta+2)x+\alpha
-\beta)y'-k(k+\alpha+ \beta+1) y; \; \; y=P_k^{( \alpha , \beta )}
(x),
\end{equation}
we obtain
\begin{equation}
\label{difeqg4} 4(1-x^2)^2  z'' = 4x(1-x^2)z' -
\end{equation}
$$
 \left[(2k+ \alpha+ \beta+1)^2
(1-x^2)-2(1+x)\alpha^2-2(1-x)\beta^2+1\right]z ; $$
$$ z=
(1-x)^{\frac{\alpha}{2}+ \frac{1}{4}}(1+x)^{\frac{\beta}{2}+
\frac{1}{4}}y, \; \; \; z^2=M_k^{\alpha }(x).
$$
Thus, in the notation of Lemma \ref{mainl},
$$A(x)=\frac{x}{2(1-x^2)} \, ,$$
$$ B(x)= \frac{(2k+ \alpha+ \beta+1)^2
(1-x^2)-2(1+x)\alpha^2-2(1-x)\beta^2+1}{4(1-x^2)^2 } \, .$$
 Now we calculate
\begin{equation}
\label{eqd}
 D=2(1-x^2)^3(4A B-B')= ( \alpha^2- \beta^2)(x^2+1)+(2 \alpha^2+2 \beta^2-1)x .
\end{equation}
\begin{theorem}
\label{thmmon} For $\alpha \ge \beta  > \frac{1}{2} \, ,$ the
consecutive maxima of the function $M_k^{\alpha , \beta }(x)$
decrease for $ x <x_0$ and increase for $ x >x_0,$ where
 $$x_0 =  \frac{\sqrt{4 \beta^2-1}-\sqrt{4 \alpha^2-1}}{\sqrt{4 \beta^2-1}+\sqrt{4 \alpha^2-1}}.$$
\end{theorem}
\begin{proof}
It is enough to show that the function $S(z;x)$ is unimodal with
the only minimum at $x_0 .$

Since $B_1=4(1-x^2) B(x),$ the numerator of $B,$ is a quadratic
with the negative leading coefficient,  by lemma \ref{mainl} it
suffices to verify that $x_0$ is the only zero of $D(x)$ in the
region defined by $B_1(x) >0.$

For, we calculate $B_1(-1)=1-4 \beta^2 \le 0,$ $B_1(1)=1-4
\alpha^2 \le 0,$
 and
$$
B_1 \left( \frac{\beta - \alpha}{\alpha+\beta+1} \right) =
\frac{(2\alpha+1)(2\beta+1)\left((2k+1)(2k+2\alpha+2\beta+1)+1
\right )}{(\alpha+\beta+1)^2}>0.
$$
Since
$$ \frac{\beta - \alpha}{\beta+\alpha+1}  \in [-1,1],$$
 $B(x)$ has precisely two zeros on $[-1,1].$

 It is easy to
check that $D$ has two real zeros for $ \alpha, \beta
> \frac{1}{2} \, , \alpha \ne \beta  . $ Moreover, for $ \alpha
\ne \beta,$
 $$ D(-1)= 1-4 \beta^2 <0, \; \; D(1)=4 \alpha^2-1 >0,$$
hence only the largest zero of $D$ lies between the zeros of
$B_1.$ If $\alpha =\beta,$ then $D=0$ implies $x=0,$ and
$$B_1
(0)=(2k+1)(2k+2\alpha+2\beta+1)+1 > 0,$$ leading to the same
conclusion. This completes the proof.
\end{proof}
\begin{remark}
Let $-1 <x_1 <...<x_k <1,$ be the zeros of $P_k^{(\alpha \beta
)}(x).$ According to Theorem \ref{thmmon} the global extremum of
$M_k^{\alpha , \beta}(x)$ lies in one of the intervals
$[\eta_{-1},x_1], \; [x_k,\eta_1],$ where $\eta_{\pm 1}$ are given
by (\ref{eqN}). Rather accurate bounds $\chi_{-1}$ and $\chi_1$ on
$x_1$ and $x_k,$ such that $x_1 < \chi_{-1} < \chi_1 <x_k,$ and
$|\eta_j- \chi_j |=O \left((k+ \alpha+ \beta)^{-2/3} \right), \;
\; j= \pm 1,$ were given in \cite{krascomp}.
\end{remark}
\section{Proof of Theorem \ref{mainth2}, even case}
\label{seceven}
 In this section we prove Theorem \ref{mainth2} for
ultraspherical polynomials of even degree.  Without loss of
generality we will assume $ x \ge 0.$

 To simplify some expressions it will be convenient to introduce the parameter
$r=2k+2\alpha+1.$

 The required differential equation  for
$$
g= (d^2-x^2)^{1/4}(1-x^2)^{\alpha/2}, \; \; g^2=M_k^{\alpha }(x,
-d , d ),
$$
 is
$$
g''-2 A(x)g'+B(x) g=0,
$$
where
$$
A(x)= \frac{x(2 d^2-1-x^2)}{2( d^2-x^2)(1-x^2)} \, ,
$$
$$
B(x)=\frac{(1-x^2)r^2-4 \alpha^2}{4(1-x^2)^2}+ \frac{2
d^2-d^4+(3-4 d^2)x^2}{4(1-x^2)( d^2-x^2)^2} \, .
$$
We also find
$$
D(x)=\frac{2( d^2-x^2)^3 (1-x^2)^2}{x} \, \left(4A B-B' \right)=$$
$$ \left( 4 \alpha^2-(1-d^2)r^2 \right)( d^2-x^2)^2+(3-4
d^2)x^4-2(5 d^4-9 d^2+3) x^2- d^6+9 d^4-9 d^2.
$$
In what follows we choose $d= \delta ,$ where $\delta$ is defined
by (\ref{defdelta}). Notice that it can be also written as
$$
\delta = \sqrt{\frac{r^2-4 \alpha^2-3}{r^2-4}} \, .
$$
The following lemma shows that $\delta$ is large enough to include
all oscillations of $M_k^{\alpha }(x).$ This fact is crucial for
our proof of Theorem \ref{mainth1}.

\begin{lemma}
\label{lmincl} The interval $(-\delta , \delta )$ contains all
local maxima of $M_k^{\alpha}(x),$ provided $\alpha >\frac{1}{2}
\, .$
\end{lemma}
\begin{proof}
The assumption $\alpha >\frac{1}{2}$ implies that $\delta$ is real
for $k \ge 0.$ It is an immediate corollary of a general result
given in  \cite{krasjac} (eq. (17) for $\lambda =0$), that in the
ultraspherical case and $k, \alpha \ge 0,$ all local maxima of
$M_k^{\alpha }(x)$ lie between the zeros of the equation
$$
A_0(x)=4k(k+2 \alpha+1)-\left((2k+2 \alpha+1)^2+4 \alpha+2
\right)x^2 =0.
$$
 Since, as easy to check, $A_0(\delta ) >0,$ the local maxima are
confined to the interval $(-\delta , \delta ).$
\end{proof}
To apply Lemma \ref{mainl} we shell check the relevant properties
of $B$ and $D,$ what will be accomplished in the following to
lemmas.
\begin{lemma}
\label{korb} Let $\alpha > \frac{1}{2} \, , \; \; k \ge 1,$ then
for $d=\delta$  the equation $B(x)=0$ has the only real positive
zero $x_0, \; \; \delta <x_0 < 1 .$ In particular, $B(x) >0$ for
$0 <x < \delta.$
\end{lemma}
\begin{proof}
It is easy to check that $r^2- 4 \alpha^2 >3, \; \; r^2 >4,$ for
$\alpha
> \frac{1}{2} \, , \; \; k \ge 1.$
The numerator $B_1$ of $B(x)$ is
 $$B_1(x)=-r^2 x^6+ \left((1+2 \delta^2)r^2+4 \delta^2-4 \alpha^2-3
 \right)x^4-$$
 $$
 \left((\delta^4+
 2 \delta^2)r^2 -\delta^4- 8 \alpha^2 \delta^2 + 6 \delta^2-3\right)x^2+
  \left( \delta^2 r^2-4 \alpha^2 \delta^2-\delta^2+2 \right)\delta^2.$$
  Using Mathematica we find the discriminant of this polynomial in
  $x,$
  $$
  Dis_x (B_1)=
  \frac{(r^2-4 \alpha^2-3) \left((r^2-4 \alpha^2-2)^2+2r^2-9 \right)(24 \alpha^2-6)^6 \,
  r^8}{(r^2-4)^{14}}\, R^2( \alpha , r),
  $$
  where
  $$R(\alpha , r)=$$
  $$ 100 (r^2-4 \alpha^2)^2 \alpha^2 r^2+7 r^6-
  (976 \alpha^2+90)r^4+(5456 \alpha^4+3180 \alpha^2+375) r^2-4(12 \alpha^2+5)^3.$$
  Under our assumptions the expressions $r^2-4 \alpha -3$ and $(r^2-4 \alpha^2-2)^2+2r^2-9
  $ are positive.
  Furthermore, rewriting $R(\alpha , r)$ in terms of $k$ and $\alpha$ one can checks that
  the substitution $\alpha \rightarrow \alpha+ \frac{1}{2}  $
  gives a polynomial consisting of monomials of the same sign.
  Thus, for any $k >0$ and $\alpha
  >\frac{1}{2}$ the discriminant does not vanish and
  the equation $B_1(x)=0$ has the same number of real zeros. For $
  \alpha=k=1$ we obtain the following test equation with just
  two real zeros,
 $$
804-2733x^2+3150x^4-1225x^6=0.
$$
It is left to demonstrate that the only positive zero $x_0$ of the
equation $B_1(x)=0,$ is in the interval $( \delta ,1 ).$ For, we
verify
$$
B_1( \delta)= 5(1-\delta^2)^2 \delta^2 >0, \; \; \; B_1( 1)= - 4
\alpha^2(1-\delta^2)^2 <0.
$$
This completes the proof.
\end{proof}
\begin{lemma}
\label{negd} Let $ \alpha > \frac{1}{2} \, , \; \; k \ge 1$ and $0
<x < \delta ,$ then
 $D(x) <0.$
\end{lemma}
\begin{proof}
We find
$$
\frac{(r^2-4)^3}{3(4 \alpha^2-1)}D(x)=2(r^2-4)(2r^2-12
\alpha^2-5)x^2 -(r^2-4 \alpha^2-3)(4r^4-4 \alpha^2-15).
$$
Then
$$
D(0) <0, \; \; D( \delta )=-5 (4 \alpha^2-1)(r^2-4 \alpha^2-3) <0,
$$
and the result follows.
\end{proof}
Applying two previous lemmas and Lemma \ref{mainl} we obtain the
following result.
\begin{lemma}
\label{lmonult}
 For $x \ge 0$ the local maxima of $M_k^{\alpha
}(x, \delta )$ form a decreasing sequence. In particular,
$\mathcal{M}_k^{\alpha }(
 \delta )=M_k^{\alpha }(0,  \delta ).$
\end{lemma}
\begin{remark}
The value of $\, \delta$ has been found as a solution of the
equation $Dis_x D =0.$ Surprisingly, it is split into linear and
biquadratic factors.  Besides trivial zeros $d=0,1,$ this equation
has four positive roots $d_1 <d_2 <d_3<d_4,$ where $d_1$ is of
order $\, O\left( \frac{1}{\sqrt{k(k+ \alpha )}} \right).$ The
other three are very close, in fact
$$d_3-d_2=O\left( \frac{1}{k^{3/2} \, \sqrt{k+ \alpha } } \right),
\; \; d_4-d_3 =O\left( \frac{\alpha^2}{k^{3/2}(k+ \alpha )^{5/2}}
\right). $$ We have chosen the simplest one $\delta =d_3.$
\end{remark}
To prove the inequality
\begin{equation}
\label{eqglavd} M_k^{\alpha } (\delta ) < \frac{2}{\pi} \,
\left(1+\frac{1}{8(k+ \alpha )^2} \right),
\end{equation}
 we have to find $M_k^{\alpha }(0,  \delta ).$
 The value of $P_k^{(\alpha , \alpha )} (0)$ for  even $k$ is (see e.g.
 \cite{abst}),
\begin{equation}
\label{ultr0}
 P_k^{(\alpha , \alpha )} (0) =
 (-1)^{k/2} \frac{\Gamma(k+\alpha+1)}{2^k \left(  \frac{k}{2} \right)!  \Gamma(\frac{k}{2}+\alpha+1) }
 \, .
\end{equation}
This yields
$$
 {\bf P}_k^{(\alpha , \alpha )} (0) = (-1)^{k/2}
 \frac{\sqrt{r \, k! \, \Gamma(r-k)}}{2^{r/2}\left(
 \frac{k}{2} \right)! \, \Gamma(\frac{r-k+1}{2}) } \, .
$$
To simplify this expression we use the following inequality (see
e.g. \cite{buis}),
\begin{equation}
\label{gamma} \frac{\Gamma(x+1)}{\Gamma^2 (\frac{x}{2}+1)}<
\frac{2^{x+\frac{1}{2}}}{\sqrt{\pi (x+ \frac{1}{2})}} \, , \; \;
\; x \ge 0,
\end{equation}
what yields for $k+ 2 \alpha \ge 0,$
$$
\left( {\bf P}_k^{(\alpha , \alpha )} (0) \right)^2 < \frac{2r}{
\pi \, \sqrt{ (2k+1)(r+2 \alpha ) }} \, .
$$

 Hence, for $|x| \le \delta,$ we have
$$\mathcal{M}_k^{\alpha }( \delta ) =M_k^{\alpha }(0, \delta )= \delta \left( {\bf P}_k^{(\alpha , \alpha )} (0)
\right)^2 < \sqrt{\frac{r^2- 4 \alpha^2}{r^2-4}} \, \frac{2r}{ \pi
\, \sqrt{ (2k+1)(r+2 \alpha ) }} . $$ It is an easy exercise to
check that for $k \ge 2, \; \alpha \ge \frac{1}{2} \, ,$ the last
expression does not exceed
$$
 \frac{2}{\pi} \,
\left(1+\frac{1}{8(k+ \alpha )^2} \right).
$$
This proves the even case of Theorem \ref{mainth2}.
\begin{remark}
In \cite{erdel} the following pointwise bound on $M_k^{\alpha ,
\beta }(x)$ is given.
\begin{equation}
\label{pointw}
 M_k^{\alpha , \beta }(x) < \frac{2e }{\pi} \;
\frac{(2k+2 \alpha+2 \beta+1) (2k+2 \alpha+2 \beta+2)}{(2k+2
\alpha+2 \beta+2)^2- \frac{2 \alpha^2 }{1-x}-\frac{2
\beta^2}{1+x}} \, .
\end{equation}
For the ultraspherical case this yields
$$M_k^{\alpha}(0)<
\frac{2e }{\pi} \left(1+O\left( \frac{\alpha^2}{k(k+ \alpha )}
\right) \right).$$ Thus, (\ref{pointw}) is quite precise, provided
$\alpha =O(k).$
\end{remark}

\section{Proof of Theorem \ref{mainth2}, odd case}
\label{secodd} In this section we will establish the odd case of
Theorem \ref{mainth2} by reducing it to the previous one. We also
give slightly more accurate bounds under the assumptions $k \ge 7,
\; \; \alpha \ge \frac{1+\sqrt{2}}{4} \, .$ They will be used in
the proof of Theorem \ref{mainth1} in the next section.

 As $\delta $ is a
function of $k$ and $\alpha ,$ to avoid ambiguities or a messy
notation arising when they vary, throughout this section we will
use  $\delta (k, \alpha )$ instead of $\delta$ and set
$\mathcal{F}_k^{\alpha} = \mathcal{M}_k^{ \alpha }( \delta),$ and
$F_k^{\alpha} (x)= M_k^{ \alpha }(x, \delta).$

 Since the value of the first, nearest to zero, maximum of
$F_k^{\alpha} (x),$ which we assume is attained at $x= \xi,$ is
unknown for odd $k$, we need some technical preparations. First of
all we have to find an upper bound on $\xi .$ Let $k=2i+1$ be odd,
and let $0=x_0 <x_1<...<x_i ,$ be the nonnegative zeros of
$P_k^{(\alpha , \alpha )} (x).$ Obviously, $0 < \xi <x_1,$ so we
can use an upper bound on $x_1 $ instead. An appropriate estimate
for zeros of ultraspherical polynomials has been given in
\cite{elbert}, in particular
$$
x_1 < \left( \frac{2k^2+1}{4k+2} + \alpha \right)^{-1/2} h_k,
$$
where $h_k$ is the least positive zero of the Hermite polynomial
$H_k (x).$

Since $h_k \le \sqrt{\frac{21}{4k+2}} \, ,$ \cite[sec.
6.3]{szego}, we obtain
\begin{equation}
\label{eqxi} \xi \le \sqrt{\frac{21}{2k^2+4 \alpha k+2\alpha
+1}}:= \xi_0.
\end{equation}
Using the formula
$$
\frac{d}{d x}P_k^{( \alpha , \beta  )}(x)= \frac{k+ \alpha +
\beta+1}{2} P_{k-1}^{( \alpha +1, \beta +1 )}(x),
$$
which for the ultraspherical orthonormal case yields
$$
\frac{d}{d x}{\bf P}_k^{( \alpha , \alpha  )}(x)= \sqrt{(r-k)k} \;
{\bf P}_{k-1}^{( \alpha+1 , \alpha +1 )}(x)
$$
and the simplest Taylor expansion around zero,
$$
{\bf P}_k^{( \alpha , \alpha  )}(\xi)=\sqrt{(r-k)k} \; {\bf
P}_{k-1}^{( \alpha+1 , \alpha +1 )}(\epsilon \xi) \, \xi , \; \;
\; 0 < \epsilon <1,
$$
what reduces the problem to the even case, we obtain
$$
F_k^{ \alpha }(\xi ) < \sqrt{\delta^2(k, \alpha )-\xi^2 } \; (1-
\xi^2 )^{\alpha } \left( {\bf P}_{k-1}^{( \alpha+1 , \alpha +1
)}(\epsilon \xi) \right)^2  (r-k)k  \, \xi^2 <
$$
$$
\frac{\sqrt{\, \delta^2(k, \alpha ) - \xi^2} \; (1-\xi^2)^{\alpha}
}{\sqrt{\, \delta^2(k-1, \alpha+1 ) - \epsilon^2 \xi^2 } \;
(1-\epsilon^2 \xi^2)^{\alpha+1}}  \; F_{k-1}^{ \alpha+1}(\epsilon
\xi)  (r-k)k \, \xi_0^2  <
$$
$$
\frac{\sqrt{ \, \delta^2(k, \alpha ) - \xi^2}}{(1- \xi^2 ) \,
\sqrt{\, \delta^2(k-1, \alpha+1 ) - \xi^2}} \; \mathcal{F}_{k-1}^{
\alpha+1} \; (r-k)k  \, \xi_0^2  \, .
$$
The last function increases in $\xi$ and substituting $\xi_0$ we
have
\begin{equation}
\label{eqtau}
 F_k^{\alpha }( \xi )< {\it{v}}(k, \alpha )
\mathcal{F}_{k-1}^{ \alpha+1},
\end{equation}
 where
$$
{\it{v}}(k, \alpha )= \frac{(r-k)k \, \xi_0^2  \; \sqrt{ \,
\delta^2(k, \alpha ) - \xi_0^2}}{(1- \xi_0^2 ) \, \sqrt{ \,
\delta^2 (k-1, \alpha+1 ) - \xi_0^2}} \, .
$$
We have checked using Mathematica that
$${\it{v}}_1 ( k, \alpha ) =\left(1+ \frac{1}{8(k+ \alpha )^2} \right){\it{v}} ( k, \alpha )$$
is a decreasing function in $k$ and $\alpha ,$ provided $k \ge 3$
and $ \alpha \ge \frac{1}{2} \, $ (an explicit expression for
$\it{v}$ is somewhat messy and is omitted). In fact, this is much
easier than one may expect as the numerator and the denominator of
$ \frac{d}{d \alpha} {\it{v}}_1^2 ( k+3, \alpha+\frac{1}{2} )$ and
$ \frac{d}{d k} {\it{v}}_1^2 ( k+3, \alpha+\frac{1}{2})$ consist
of the monomials of the same sign.

 Calculations yield
  $$ {\it{v}}_1( 3,
\frac{1}{2})< 115 \, , \; \; \; {\it{v} }_1(7,\frac{1+
\sqrt{2}}{4} ) < \frac{29}{2} \, .$$

Finally, applying (\ref{eqglavd}) and (\ref{eqtau}) and coming
back to the usual notation, we conclude
\begin{lemma}
Let $k$ be odd, then
\begin{equation}
\label{odl} \mathcal{M}_{k}^{ \alpha}( \delta) \le \left\{
\begin{array}{ccc}
 \frac{230}{\pi} \, &
k \ge 3, & \alpha > \frac{1}{2} \,  ,\\
& & \\
\frac{29}{\pi} \, ,& k \ge 7, & \alpha > \frac{1+ \sqrt{2}}{4} \,
.
\end{array}
\right.
\end{equation}
\end{lemma}
  This completes the proof of Theorem
\ref{mainth2}.

\section{Proof of Theorem \ref{mainth1}}
\label{secth1}
 First, we will establish the following bounds which
are slightly better than these of Theorem \ref{mainth1} but stated
in terms of $r=2k+2 \alpha+1,$ and $\tau = \frac{2 \alpha +1}{r}
\, .$ It is worth noticing that in some respects $r$ and $\tau$
are more natural parameters than $k$ and $\alpha $ (see
\cite{krasjac}).
\begin{lemma}
\label{glav}
\begin{equation}
\label{eqglav} \mathcal{M}_k^{\alpha} < \left\{
\begin{array}{cc}
\frac{12}{13} \, r^{1/3} \tan^{1/3} \tau , & k \ge 6, \; \; even,\\
& \\
14 \, r^{1/3} \tan^{1/3} \tau , & k \ge 7, \; \; odd.
\end{array}
\right.
\end{equation}
provided $k \ge 6, \; \; \alpha \ge \frac{1+ \sqrt{2} }{4} \, .$
\end{lemma}
\begin{proof}
 Let $\epsilon = \frac{
2^{-1/3}}{3} \,r^{-2/3} \, tan^{4/3} \tau .$ It is easy to check
that  $ \epsilon < \frac{1}{31} \, ,$ (the extremal case
corresponds to $k=6, \; \alpha = \infty $).

Since
$$\delta > \cos  \tau > \eta = \left(1 -  \epsilon
\right)\cos \tau ,$$ where $\eta$ is defined in (\ref{eqNus}), it
follows by Theorem \ref{grmax} that all local maxima of
$M_k^{\alpha} (x )$ are inside the interval $(- \delta, \delta) .$
Now we have
\begin{equation}
\label{eqmax}
 \max_{|x| \le 1} \left\{(1-x^2)^{\alpha +
\frac{1}{2}}\left({\bf P}_k^{(\alpha, \alpha )} (x) \right)^2
\right\}= \mathcal{M}_k^{\alpha} (\delta ) \max_{0 \le x \le \eta}
\sqrt{\frac{1-x^2}{\delta^2-x^2}} =
\end{equation}
$$
\mathcal{M}_k^{\alpha} (\delta ) \,
\sqrt{\frac{1-\eta^2}{\delta^2-\eta^2}} \, .
$$
 By the explicit expression for $\epsilon$ given by (\ref{eqNus}), one can check
 that the function $\sqrt{ \, 2- \epsilon} \,$ increases in $k$ and a decreases in $\alpha .$
 We obtain by $\epsilon < \frac{1}{31} \, ,$
$$
\sqrt{\delta^2-\eta^2} > \sqrt{\cos^2 \tau-\eta^2} =
\sqrt{\epsilon (2- \epsilon ) } \, \cos \tau > \frac{7}{5}
\,\sqrt{\epsilon } \, \cos \tau  .
$$
 Using the restrictions $k \ge 6, \; \alpha \ge \frac{1+\sqrt{2}}{4} \,
 ,$ and a simple trigonometric inequality,
 we find
$$
\sqrt{1-\eta^2}=\sqrt{1-(1- \epsilon )^2 \cos^2 \tau } \le \sin
\tau \left( 1+ \epsilon \, \cot^2 \tau  \right) =$$
$$
\left( 1+  \frac{1}{3} \, \left(\frac{2k(k+ 2 \alpha +1)}{(2
\alpha+1 )^2 (2k+ 2 \alpha +1)^2} \right)^{1/3}\right)\sin \tau <
\frac{37}{32} \, \sin \tau .
$$
 Thus, we obtain
$$
\sqrt{\frac{1-\eta^2}{\delta^2-\eta^2}}< \frac{185 \tan \tau}{224
\sqrt{ \epsilon} }= \frac{185  \sqrt{3} }{224  } r^{1/3}
\tan^{1/3} \tau  < \frac{13}{9} \, r^{1/3} \tan^{1/3} \tau ,$$ and
the result follows by (\ref{eqmax}) and (\ref{eqglavd}) for $k$
even, and (\ref{odl}) for $k$ odd.
\end{proof}
Now Theorem \ref{mainth1} is an immediate corollary of
(\ref{eqglav}) and
$$
\frac{r^{1/3} \tan^{1/3} \tau }{ \alpha^{1/3} \left(1+ \frac{
\alpha }{k} \right)^{1/6}} = \left( \frac{(2 \alpha+1)^2 (2k+2
\alpha+1 )^2}{4 \alpha^2 (k+ \alpha)(k+2 \alpha+1) }\right)^{1/6}
\le ( 4 \sqrt{2}-2 )^{1/3},
$$
for $\alpha \ge \frac{1+ \sqrt{2}}{4} \, .$ This completes the
proof.

\hspace{1pt}

 \noindent
 {\bf Acknowledgement.} I am grateful to D.K. Dimitrov for a
 helpful discussion and, especially, for bringing my attention to Sonin's function.


\end{document}